\documentclass{amsart}
\usepackage{amsfonts,amssymb}
\usepackage{graphicx}
\usepackage[noadjust]{cite}

 4

\usepackage[british]{babel}
\usepackage{amsmath,bbm}

\usepackage{kpfonts}  

\usepackage[hidelinks]{hyperref}
\usepackage[T1]{fontenc}
\usepackage[cp1250]{inputenc}
\usepackage{float}
\usepackage{enumitem}
\usepackage{mathabx}
\usepackage[activate={true,nocompatibility},final,tracking=true,kerning=true,spacing=true,factor=1100,stretch=10,shrink=10]{microtype}
\usepackage{upgreek}

\newcommand{\assign}{:=}
\newcommand{\nocomma}{}
\newcommand{\tmem}[1]{{\em #1\/}}

\newcommand{\tmop}[1]{\ensuremath{\operatorname{#1}}}

\newcommand{\tmtextit}[1]{{\itshape{#1}}}

\newtheorem{lem}{Lemma}[section]
\newtheorem{cor}[lem]{Corollary}

\newtheorem{thm}[lem]{Theorem}
\newtheorem{prop}[lem]{Proposition}

\newtheorem{rk}{Remark}
\newtheorem*{problem**}{Problem}
\newtheorem*{thm**}{Theorem}

\newtheorem*{ex**}{Example}
\newtheorem*{rk**}{Remark}

\newenvironment{prf}{\noindent\textbf{Proof.\ }}{\hspace*{\fill}$\Box$\medskip}

\title{Non-degenerate jump of Milnor numbers of surface singularities}
\author{Szymon Brzostowski, Tadeusz Krasi{\'n}ski and Justyna Walewska}
\keywords{Milnor number, deformation of singularity, non-degenerate singularity, Newton polyhedron}
\subjclass[2010]{14B07, 32S30}
\begin{document}
\maketitle

\renewcommand{\thepage}{[\arabic{page}]}

\begin{abstract}
The jump of the Milnor number of an isolated singularity $f_0$ is the minimal non-zero difference between the Milnor numbers of $f_0$ and one of its deformations $(f_s)$. We give a formula for the jump in some class of surface singularities in the case deformations are non-degenerate.
\end{abstract}

\vspace{3mm}

\section{Introduction}\label{Par 1}

Let $f_0 : (\mathbbm{C}^n, 0) \rightarrow (\mathbbm{C}, 0)$ be an
{\tmem{(isolated) singularity}}, i.e.~let $f_0$ be a germ at $0$ of a
holomorphic function having an isolated critical point at $0 \in
\mathbbm{C}^n$, and $0 \in \mathbbm{C}$ as the corresponding critical value.
More specifically, there exists a representative $\hat{f}_0 : U \rightarrow
\mathbbm{C}$ of $f_0$ holomorphic in an open neighborhood $U$ of the point $0
\in \mathbbm{C}^n$ such that:

\begin{itemize}
  \item[$\bullet$] $\hat{f}_0 (0) = 0$,
  
  \item[$\bullet$] $\nabla \hat{f}_0 (0) = 0$,
  
  \item[$\bullet$] $\nabla \hat{f}_0 (z) \neq 0$ for $z \in U\backslash \{ 0 \}$,
\end{itemize}

{\noindent}where for a holomorphic function $f$ we put $\nabla f \assign
(\partial f / \partial z_1, \ldots, \partial f / \partial z_n)$.

In the sequel we will identify germs of functions with their representatives
or the corresponding convergent power series. The ring of germs of holomorphic functions of $n$ variables will be denoted by $\mathcal{O}_n$.

A {\tmem{deformation of the singularity}} $f_0$ is a germ of a holomorphic
function $f = f (s, z) : (\mathbbm{C} \times \mathbbm{C}^n, 0) \rightarrow
(\mathbbm{C}, 0)$ such that:

  \begin{itemize}
  \item[$\bullet$] $f (0, z) = f_0(z)$,
  
  \item[$\bullet$] $f (s, 0) = 0$,

\end{itemize}
The deformation $f (s, z)$ of the singularity $f_0$ will also be treated as a
family $(f_s)$ of germs, putting $f_s (z) \assign f (s, z)$. Since $f_0$ is an isolated singularity, $f_s$ has also isolated singularities near the origin, for sufficiently small $s$ {\cite[Theorem 2.6 in Chap.~I]{GLS07}}.
\begin{rk**}
  \label{nota1}Notice that in the deformation $(f_s)$ there can occur in
  particular {\tmem{smooth}} germs, that is germs satisfying $\nabla f_s (0)
  \neq 0$. In this context,
the symbol $\nabla f_s$ will always denote $\nabla_z f_s(z)$.

\end{rk**}

By the above assumptions it follows that, for every sufficiently small $s$,
one can define a (finite) number ${\upmu}_s$ as the Milnor number of $f_s$,
namely
\[ {\upmu}_s \assign {\upmu} (f_s) \mathrel{=} \dim_{\mathbbm{C}} 
   \mathbin{\mathbin{\mathcal{O}}_n / \mathbin{(\nabla f_s)}}
   \mathord{\mathrel{}} \mathrel{=} {\upmu} \left( \frac{\partial
   f}{\partial z_1}, \ldots, \frac{\partial f}{\partial z_n} \right) \text{,}
\]
where the symbol ${\upmu} ( \frac{\partial f}{\partial z_1}, \ldots,
\frac{\partial f}{\partial z_n} )$ denotes intersection multiplicity of
the ideal $( \frac{\partial f}{\partial z_1},\allowbreak \ldots, \frac{\partial
f}{\partial z_n}) \mathcal{O}_n$ in $\mathcal{O}_n$.

Since the Milnor number is upper semi-continuous in the Zariski topology in families of singularities
{\cite[Theorem 2.57 in Chap.~II]{GLS07}}, there exists an open neighborhood
$S$ of the point $0 \in \mathbbm{C}$ such that

 \begin{itemize}
  \item[$\bullet$] ${\upmu}_s = \tmop{const.}$ for $s \in S \setminus \{ 0 \}$,
  
  \item[$\bullet$] ${\upmu}_0 \geqslant {\upmu}_s$ for $s \in S$.
\end{itemize}

The (constant) difference ${\upmu}_0 -{\upmu}_s$ for $s \in S \setminus \{
0 \}$ will be called {\tmem{the jump of the deformation $(f_s)$}} and denoted
by $\uplambda ((f_s))$. The smallest nonzero value among all the jumps of
deformations of the singularity $f_0$ (such a value exists because one can
always consider a deformation of $f_0$ built of smooth germs and then for it
it is ${\upmu}_s = 0$; cf. Remark \ref{nota1}) will be called {\tmem{the jump
(of the Milnor number) of the singularity $f_0$}} and denoted by $\uplambda
(f_0)$.

The first general result concerning the jump was S.~Gusein-Zade's
{\cite{Gus93}}, who proved that there exist singularities $f_0$ for which
$\uplambda (f_0) > 1$ and that for irreducible plane curve singularities it
holds $\uplambda (f_0) = 1$. In {\cite{BK14}} the authors proved that $\uplambda
(f_0)$ is not a topological invariant of $f_0$ but it is an invariant of the
stable equivalence. The computation of $\uplambda (f_0)$ for a specific
reducible singularity (or for a class of reducible singularities) is not an
easy task. It is related to the problem of adjacency of classes of singularities. Only for a few classess of singularities we know the exact value of
$\uplambda (f_0)$. For plane curve singularities ($n=2$) we have (see {\cite{AGV85}} for terminology):
  \begin{itemize}
  \item[$\bullet$] for the one-modal family of singularities in the $X_9$ class, that is
  singularities of the form
  \[ f_0^a (x, y): = x^4 + y^4 + ax^2 y^2  \text{, {\hspace{1em}}} a\in\mathbbm{C} \text{, {\hspace{1em}}} a^2 \neq 4,
  \]
  we have $\uplambda (f_0^a) = 2$ ({\cite{BK14}}),
  
  \item[$\bullet$] for the two-modal family of singularities in the $W_{1, 0}$ class,
  that is singularities of the form
  \[ f_0^{a, b} (x, y): = x^4 + y^6 + (a + by) x^2 y^3 \text{, {\hspace{1em}}} a,b\in\mathbbm{C} \text{, {\hspace{1em}}}
     a^2 \neq 4 \text{,} \]
  we have
  \[ \uplambda (f_0^{a, b}) = \left\{ \begin{array}{ll}
       1, & \text{if } a = 0 \text{ ({\cite{BK14}})}\\
       \geqslant 2, & \text{for generic } a, b \text{ ({\cite{Gus93}})},
     \end{array} \right. \]
  \item[$\bullet$] for specific homogenous singularities $f_0^d (x, y): = x^d + y^d$, $d \geqslant
  2$, we have \newline $\uplambda (f_0^d) = \left[ \frac{d}{2} \right]$ ({\cite{BKW14}}),
  \item[$\bullet$] for homogeneous singularities of degree $d$ with generic coefficients $f_0$ we have $\uplambda(f_0)<\left[ \frac{d}{2} \right]$ ({\cite{BKW14}})
\end{itemize}

In the present paper we consider a weaker problem: {\tmem{compute the jump
$\uplambda^{\tmop{nd}} (f_0)$ of $f_0$ over all non-degenerate deformations of
$f_0$}} (i.e.~the $f_s$ in the deformations $(f_s)$ of $f_0$ are
non-degenerate singularities). Clearly, we always have $\uplambda (f_0)
\leqslant \uplambda^{\tmop{nd}} (f_0)$. Up to now, this problem has been studied
only for plane curve singularities

  \begin{itemize}
  \item[$\bullet$] A.~Bodin ({\cite{Bod07}}) gave a formula for $\uplambda^{\tmop{nd}}
  (f_0)$ for $f_0$ convenient with its Newton polygon reduced to one segment,
  
  \item[$\bullet$] J.~Walewska in \cite{Wal13} generalized Bodin's results to the non-convenient
  case,
  
  \item[$\bullet$] the authors ({\cite{BKW14}}) calculated all possible Milnor numbers of
  all non-degenerate deformations of homogenous singularities,
  
  \item[$\bullet$] J.~Walewska ({\cite{Wal10}}) proved that the {\tmem{second
  non-degenerate jump of}} $f_0$ satisfying Bodin's assumptions is equal to
  $1$.
\end{itemize}

\begin{figure}[H]
\centering
  \resizebox{150px}{100px}{\includegraphics{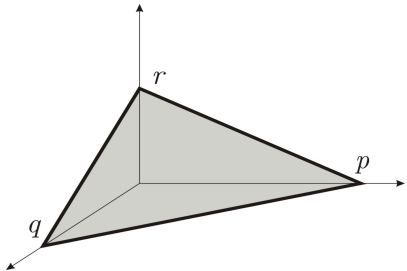}}\\
  \caption{The Newton diagram of $ f_0 (x, y, z) = ax^p + by^q + cz^r + \ldots$}\label{rys:1}
\end{figure}

In this paper we want to pass to surface singularities ($n = 3$). We give a
formula (more precisely: a simple algorithm) for $\uplambda^{\tmop{nd}} (f_0)$
in the case where $f_0$ is non-degenerate, convenient and has its Newton
diagram reduced to one triangle, (see Figure \ref{rys:1}) i.e.~$f_0$ of the form
\[ f_0 (x, y, z) = ax^p + by^q + cz^r + \ldots \hspace{1em} (p, q, r \geqslant
   2 \nocomma,\hspace{0.5em} abc \neq 0) . \]
   
\smallskip

{\noindent}Moreover, for simplicity reasons, we will only consider the case of
$p$, $q$, $r$ being pairwise coprime integers. The general case of arbitrary
$p$, $q$, $r$ will be the topic of a next paper.\label{Chyba raczej nie}

\section{Non-degenerate singularities}

In this Section we recall the notion of non-degenerate singularities. We restrict ourselves to surface singularities. All notions can easily be generalized to higher dimensions.
\noindent Let $f_0 (x, y, z) \assign \sum_{i, j, k \in \mathbbm{N}} a_{i j k} x^i y^j
z^k$, be a singularity. Let $\tmop{supp} (f_0) \assign \{
(i, j, k) \in \mathbbm{N}^3 : a_{i j k} \neq 0 \}$ be the {\tmem{support of
$f_0$}}. The {\tmem{Newton polyhedron $\Gamma_+ (f_0)$ of $f_0$}} is the
convex hull of the set
\[ \bigcup_{(i, j, k) \in \tmop{supp} (f_0)} (i, j, k) +\mathbbm{R}_+^3, \]
where $\mathbbm{R}^3_+$ is the closed octant of $\mathbbm{R}^3$ consisting of
points with nonnegative coordinates. The boundary (in $\mathbbm{R}^3$) of
$\Gamma_+ (f_0)$ is an unbounded polyhedron with a finite number of
$2$-dimensional faces, which are (not necessarily compact) polygons. The
singularity $f_0$ is called {\tmem{convenient}} if $\Gamma_+ (f_0)$ has some
points in common with all three coordinate axes in $\mathbbm{R}^3$. The set of
compact faces (of all dimensions) of $\Gamma_+ (f_0)$ constitutes the
{\tmem{Newton diagram of $f_0$}} and is denoted by $\Gamma (f_0)$. For each
face $S \in \Gamma (f_0)$ we define a weighted homogeneous polynomial
\[ (f_0)_S \assign \sum_{(i, j, k) \in S} a_{i j k} x^i y^j z^k . \]
We call the singularity $f_0$ {\tmem{non-degenerate on $S \in \Gamma (f_0)$}}
if the system of equations
\[ \frac{\partial (f_0)_S}{\partial x} (x, y, z) = 0, \hspace{1em}
   \frac{\partial (f_0)_S}{\partial y} (x, y, z) = 0, \hspace{1em}
   \frac{\partial (f_0)_S}{\partial z} (x, y, z) = 0 \]
has no solutions in $(\mathbbm{C}^{\ast})^3$; $f_0$ is {\tmem{non-degenerate}}
({\tmem{in the Kouchnirenko sense}}) if $f_0$ is non-degenerate on every face
$S \in \Gamma (f_0)$.

\begin{figure}[H]
\begin{center}
  
  \resizebox{362px}{76.5px}{\includegraphics{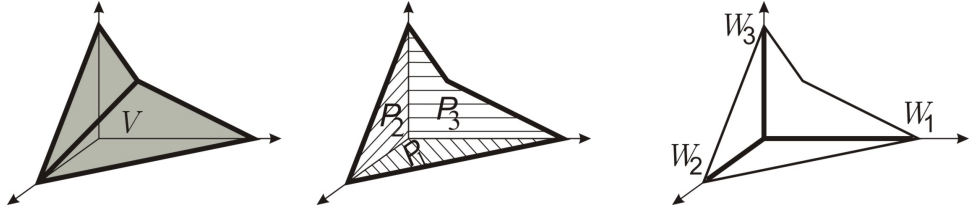}}
  \caption{Geometric meaning of volume $V$, areas $P_i$ and lengths $W_j$.}
  \label{rys:2}
\end{center}
\end{figure}

Assume now that $f_0$ is convenient. We introduce the following notation:
\begin{itemize}
  \item[$\bullet$] $\Gamma_- (f_0)$ -- the compact polyhedron bounded by $\Gamma (f_0)$ and the three coordinate planes (labeled in a self-explanatory way as
  $\tmop{OXY}$, $\tmop{OXZ}$, $\tmop{OYZ}$); in other words, $\Gamma_- (f_0)
  \assign \overline{\mathbbm{R}_+^3 \setminus \Gamma_+ (f_0)}$,
  
   \item[$\bullet$] $V$ -- the volume of $\Gamma_- (f_0)$,
  
   \item[$\bullet$] $P_1$, $P_2$, $P_3$ -- the areas of the two-dimensional faces of
  $\Gamma_- (f_0)$ lying in the planes $\tmop{OXY}$, $\tmop{OXZ}$,
  $\tmop{OYZ}$, respectively; e.g.~$P_1$ is the area of the set $\Gamma_-
  (f_0) \cap \tmop{OXY}$,
  
   \item[$\bullet$] $W_1$, $W_2$, $W_3$ -- the lengths of the edges (= one-dimensional
  faces) of $\Gamma_- (f_0)$ lying in the axes $\tmop{OX}$, $\tmop{OY}$,
  $\tmop{OZ}$, respectively (see Figure \ref{rys:2} ).
\end{itemize}

We define the {\tmem{Newton number $\upnu (f_0)$}} of $f_0$ by
\[ \upnu (f_0) \assign 3! V - 2! (P_1 + P_2 + P_3) + 1! (W_1 + W_2 + W_3) - 1.\tag{$\circ$}\label{o}
\]
The importance of $\upnu (f_0)$ has its source in the celebrated {\tmem{Kouchnirenko
theorem}}:

\begin{thm}[\cite{Kou76}]If $f_0$ is a convenient
singularity, then
\begin{enumerate}
  \item ${\upmu} (f_0) \geqslant \upnu (f_0)$,
  
  \item if $f_0$ is non-degenerate then ${\upmu} (f_0) = \upnu (f_0)$.
\end{enumerate}
\end{thm}

\begin{rk}The Kouchnirenko theorem is true in any dimension {\cite{Kou76}}. 
\end{rk}

\section{Non-degenerate jump of Milnor numbers of singularities}
Let $f_0\in\mathcal{O}_3$ be a singularity. A deformation $(f_s)$ of $f_0$ is called \textit{non-degenerate} if $f_s$ is non-degenerate for  $s\neq 0$. The set of all non-degenerate deformations of the singularity $f_0$ will be denoted by $\mathcal{D}^{\operatorname{nd}}(f_0)$. \textit{Non-degenerate jump} $\uplambda^{\operatorname{nd}}(f_0)$ \textit{of the singularity} $f_0$ is the minimal of non-zero jumps over all non-degenerate deformations of $f_0$, which means
$$\uplambda^{\operatorname{nd}}(f_0):=\min_{(f_s)\in \mathcal{D}^{\operatorname{nd}}_0 (f_0)} \uplambda ((f_s)),$$
where by $\mathcal{D}^{\operatorname{nd}}_0 (f_0)$ we denote all the non-degenerate deformations $(f_s)$ of $f_0$ for which $\uplambda ((f_s))\neq 0$.

Obviously
\begin{prop}
For each singularity $f_0$ we have the inequality 
$$\uplambda (f_0)\leq \uplambda^{\operatorname{nd}}(f_0).$$
\end{prop}
In investigations concerning $\uplambda^{\operatorname{nd}}(f_0)$ we may restrict our attention to non-degenerate $f_0$ because the non-degenerate jump for degenerate singularities can be found using the proposition below (cf.~\cite[Lemma 5]{Bod07}).
Let $f^{\operatorname{nd}}_0$ denote any non-degenerate singularity for which $\Gamma(f_0)=\Gamma(f^{\operatorname{nd}}_0)$. Such singularities always exist.

\begin{prop}\label{wzorklam}
If $f_0$ is degenerate then

$$\uplambda^{\operatorname{nd}}(f_0)=\bigg\{
\begin{array}{ll}
\upmu(f_0)-\upmu(f^{\operatorname{nd}}_0), & \text{if }\upmu(f_0)-\upmu(f^{\operatorname{nd}}_0)>0\\
\uplambda^{\operatorname{nd}}(f^{\operatorname{nd}}_0), & \text{if }\upmu(f_0)-\upmu(f^{\operatorname{nd}}_0)=0
\end{array}.$$
\end{prop}
\begin{prf}
 This follows from the fact that a generic small perturbation of coefficients of these monomials of $f_0$ which correspond to points belonging to $\bigcup\Gamma(f_0)$ (which are finite in number) give us non-degenerate singularities with the same Newton polyhedron as $f_0$.
\end{prf}

\begin{rk}
By the P\l oski theorem (\cite[Lemma 2.2]{Plo90}, \cite[Theorem 1.1]{Plo99}), for degenerate plane curve singularities ($n=2$) the second possibility in Proposition \ref{wzorklam} is excluded.
\end{rk}

A crucial r\^ole in the search for the formula for $\uplambda^{\operatorname{nd}}(f_0)$ will be played by the monotonicity of the Newton number with respect to the Newton polyhedron. Namely, J.~Gwo\'zdziewicz \cite{Gwo08} and M.~Furuya \cite{Fur04} proved:

\begin{thm}[Monotonicity Theorem]
Let $f_0, \tilde{f_0}\in\mathcal{O}_n$ be two convenient singularities such that $\Gamma_{+}(f_0)\subset \Gamma_{+}(\tilde{f_0})$. Then $\upnu(f_0)\geqslant \upnu(\tilde{f_0})$.
\end{thm}
 
By this theorem the problem of calculation of $\uplambda ^{\operatorname{nd}}(f_0)$ can be reduced to a purely combinatorial one. Namely, we define specific deformations of a convenient and non-degenerate singularity $f_0\in\mathcal{O}_n$. Denote by $J$ the set of integer points $\textit{\textbf{i}}=(i_1,\ldots,i_n)\neq 0$ lying in the closed domain bounded by coordinate hyperplanes $\{z_i=0\}$ and the Newton diagram; in other words $J:=\Gamma_- (f_0)\cap\mathbbm{Z}^n$. Obviously, $J$ is a  finite set. For $\textit{\textbf{i}}=(i_1,\ldots,i_n)\in J$ we define the deformation $(f^{\textit{\textbf{i}}}_s)_{s\in\mathbbm{C}}$ of $f_0$ by the formula

$$f^{\textit{\textbf{i}}}_s(z_1,\ldots,z_n):=f_0(z_1,\ldots,z_n)+sz^{i_1}_1\ldots z^{i_n}_n.$$

\begin{prop}
For every $\textit{\textbf{i}}\in J$ the deformation $(f^{\textit{\textbf{i}}}_s)$ of $f_0$ is convenient and non-degenerate for all sufficiently small $|s|$.
\end{prop}

\begin{prf}
See \cite{Kou76} or \cite[Appendix]{Oka79}.

\end{prf}

Combining the Monotonicity Theorem with the above proposition we reach the conclusion that in order to find $\uplambda^{\operatorname{nd}}(f_0)$ it is enough to consider only the non-degenerate deformations of the type $(f^{\textit{\textbf{i}}}_s)$.
\begin{thm}\label{twieA}
If $f_0$ is a convenient and non-degenerate singularity, then
$$\uplambda^{\operatorname{nd}}(f_0)=\min_{\textit{\textbf{i}}\in J_0}\uplambda((f^{\textit{\textbf{i}}}_s))$$
where $J_0\subset J$ is the set of these $\textit{\textbf{i}}\in J$ for which $\uplambda^{\operatorname{nd}} ((f^{\textit{\textbf{i}}}_s))>0$.
\end{thm}

\begin{prf}
By the Kouchnirenko theorem it suffices to consider non-degenerate deformations of $f_0$ of the form

\begin{equation}\tag{$\Asterisk$}\label{wzor*}
f_s(z_1,\ldots,z_n)=f_0(z_1,\ldots,z_n) + \sum_{\textit{\textbf{i}}\in J}a_{\textit{\textbf{i}}}(s)z^{\textit{\textbf{i}}},
\end{equation}
where $a_{\textit{\textbf{i}}}(s)$ are holomorphic at $0\in\mathbbm{C}$ and $a_{\textit{\textbf{i}}}(0)=0$. Then by the Monotonicity Theorem we may restrict the scope of deformations (\ref{wzor*}) to deformations with only one term added i.e.~the deformations $(f^{\textit{\textbf{i}}}_s)$ for $\textit{\textbf{i}}\in J_0$.
\end{prf}

\begin{cor}\label{wnio1}
If $f_0$ and $\tilde{f_0}$ are non-degenerate and convenient singularities and $\Gamma(f_0)=\Gamma(\tilde{f_0})$ then $\uplambda^{\operatorname{nd}}(f_0)=\uplambda^{\operatorname{nd}}(\tilde{f_0})$.
\end{cor}

\section{An algorithm for $\uplambda^{\operatorname{nd}}(f_0)$ in the case of one face Newton diagram of surface singularities}

In this Section we give a simple algorithm for calculating $\uplambda^{\operatorname{nd}}(f_0)$ provided that $f_0\in\mathcal{O}_3$ is a convenient and non-degenerate singularity with one two-dimensional face of its Newton diagram. Let $p,q,r$ be the first (i.e.~nearest to the origin) points of $\Gamma_{+}(f_0)$ lying on the axes $\tmop{OX}, \tmop{OY}$ and $\tmop{OZ}$, respectively. Then by Corollary \ref{wnio1} we may assume that
$$f_0(x,y,z)=x^p+y^q+z^r,\qquad p,q,r\geqslant 2.$$
By formula (\ref{o}) we have $\upmu(f_0)=(p-1)(q-1)(r-1)$. Moreover, without loss of generality we may also assume that 

\begin{equation}\tag{$\dagger$}\label{rown+}
p\geqslant q\geqslant r.
\end{equation}
Additionally, we demand that $p,q,r$ are pairwise coprime

\begin{equation}\tag{$\Asterisk\Asterisk$}\label{wzor**}
\operatorname{GCD}(p,q)=\operatorname{GCD}(p,r)=\operatorname{GCD}(q,r)=1.
\end{equation}

By Theorem \ref{twieA} we have to compare the jumps of deformations $(f^{\textit{\textbf{i}}}_s)_{s\in\mathbbm{C}}$, where $\textit{\textbf{i}}\in J$, i.e.~$\textit{\textbf{i}}$ are integer points lying in the octant of $\mathbbm{R}^3$ under the triangle with vertices $(p,0,0)$, $(0,q,0)$, $(0,0,r)$ (see Figure \ref{rys:1}).

\begin{enumerate}
\item[I.] First we consider points in $J$ lying on the axes. Using formula (\ref{o}) and assumption (\ref{rown+}) we easily check that the axes-jump is realized by the deformation $(f^{(p-1,0,0)}_s)$, i.e.
$$f^{(p-1,0,0)}_s(x,y,z)=x^p+y^q+z^r+sx^{p-1},$$
and the jump is equal to $(q-1)(r-1)$.

\item[II.] Now we consider points in $J$ lying in coordinate planes. By the results of Bodin \cite{Bod07} and Walewska \cite{Wal10} we easily check that the minimal jumps on respective planes are realized by

\begin{enumerate}
\item [i.] the deformation $(f^{(b_1,q-a_1,0)})$, where $a_1,b_1\in\mathbbm{Z}$ are such that $a_1p-b_1q=1$ and $0<a_1<q$, $b_1>0$; this delivers the $OXY$-jump equal to $(r-1)$,

\item [ii.] the deformation $(f^{(0,b_2,r-a_2)}_s)$, where $a_2,b_2\in\mathbbm{Z}$ are such that $a_2q-b_2r=1$ and $0<a_2<r$, $b_2>0$; this delivers the $OYZ$-jump equal to $(p-1)$,

\item [iii.] the deformation $(f^{(b_3,0,p-a_3)}_s)$, where $a_3,b_3\in\mathbbm{Z}$ are such that $a_3p-b_3r=1$ and $0<a_3<p$, $b_3>0$; this delivers the $OXZ$-jump equal to $(q-1)$.
\end{enumerate}
\end{enumerate}

The above considerations imply that the jump realized by the points lying either in coordinate planes or on axes is equal to $(r-1)$.
\begin{enumerate}
\item [III.] Let us pass to the deformations $(f^{\textit{\textbf{i}}}_s)$ for which the point $\textit{\textbf{i}}$ lies in the interior of the tetrahedron with vertices $(0,0,0)$, $(p,0,0)$, $(0,q,0)$, $(0,0,r)$. Any such point $(\alpha, \beta, \gamma)$ satisfies the conditions:

\begin{enumerate}[label=(\Alph*)]
\item\label{war1} $0<\alpha<p$, $0<\beta<q$, $0<\gamma<r$,
\item\label{war2} $\dfrac{\alpha}{p}+\dfrac{\beta}{q}+\dfrac{\gamma}{r}<1 \text{ }\text{or equivalently} \text{ }\alpha qr+\beta pr+\gamma pq<pqr$.
\end{enumerate}
Moreover, the jump of the deformation $(f^{(\alpha,\beta,\gamma)}_s)$ is equal to $6$ times the volume of the tetrahedron with vertices $(p,0,0)$, $(0,q,0)$, $(0,0,r)$, $(\alpha,\beta,\gamma)$ i.e.~
$$pqr-\alpha qr - \beta pr - \gamma pq.$$
\end{enumerate}

Thus, we have reduced our original problem to a number theoretic one.

\begin{problem**}
Given pairwise coprime integers $p>q>r$ greater than $1$. Find positive integers $\alpha$, $\beta$, $\gamma$ satisfying \ref{war1} and \ref{war2} for which the expression $pqr-\alpha qr - \beta pr - \gamma pq$ attains its positive minimum.
\end{problem**}

In order to solve it, first notice that $\operatorname{GCD}(qr,pr,pq)=1$. Consequently, there are integers $a$, $b$, $c$ such that
\begin{equation}\tag{$\ddagger$}\label{rown++}
aqr+bpr+cpq=1.
\end{equation}
They can be obtained by the Euclid algorithm using the well-known associativity law: \tmtextit{for any integers $u$, $v$, $w$ we have $\operatorname{GCD}(u,v,w)=\operatorname{GCD}(\operatorname{GCD}(u,v),w)$}. Notice that in any identity of the type (\ref{rown++}) it holds $abc\neq 0$. If we write $a=a'p+a''$, $0\leqslant a''<p$, then, by abuse of notation, we obtain yet another identity $aqr + bpr + cpq = 1$, but now $0<a<p$.
Next, we write $b=b'q-b''$, $0< b''<q$, and we use it to obtain a similar identity $aqr - bpr + cpq = 1$ in which $0<a<p$ and $0<b<q$.
Notice that then $0<|c|<r$. In fact, $|cpq|=|1-aqr+bpr|\leqslant 1+ r|bp-aq|\leqslant 1 + r(pq-p-q)=pqr-pr-qr+1 < pqr$. Thus, finally we have obtained the identity

\begin{equation}\tag{$\square$}\label{rownkwadrat}
aqr-bpr+cpq=1,  \text{ where } 0<a<p, \text{ }  0<b<q,  \text{ } 0<|c|<r.
\end{equation}

Now we consider two cases:

\begin{enumerate}
  \item $c < 0$. Then the triple $\alpha=p-a$, $\beta=b$, $\gamma=-c$ is the solution that we seek for.
  In fact, $\alpha, \beta, \gamma$  clearly satisfy \ref{war1}, moreover $pqr-\alpha qr-\beta pr-\gamma pq=aqr-bpr+cpq=1$.
  This is the optimal value one can hope for, so the Problem is solved in this case. Hence $\uplambda^{\tmop{nd}}(f_0)=1$ and the deformation $(f^{p-a,b,-c}_s)$ realizes the jump $1$.
  \item \label{c-less-0}$c > 0$. Under this condition, we claim that there is no point
  $(\alpha, \beta, \gamma)$ satisfying both \ref{war1} and \ref{war2} and for
  which the minimum in the Problem is equal to $1$. In fact, if there existed
  such a point, then from the relation $pqr - \alpha qr - \beta pr - \gamma pq
  = 1$ we would get $(p - \alpha) qr - \beta pr - \gamma pq = 1$, which together with
  (\ref{rownkwadrat}) would imply that $(p - (\alpha + a)) qr = (\beta - b) pr + (\gamma
  + c) pq$. But since $\tmop{GCD} (p, r) = \tmop{GCD} (p, q) = 1$ and $| p -
  (\alpha + a) | < p$, this is only possible when $\alpha = p - a$. Hence, we
  would get $(\beta - b) r + (\gamma + c) q = 0$. Similarly, since $\tmop{GCD}
  (r, q) = 1$ and $| \beta - b | < q$, we would obtain $\beta = b$ and
  consequently $\gamma = - c < 0$, contradictory to \ref{war1}.
\end{enumerate}
The above observation means that in case (\ref{c-less-0}) we must further
continue our search for $\alpha, \beta, \gamma$ solving the Problem.
Accordingly, we repeat the above reasoning for the identity
\[ aqr + bpr + cpq = 2, \]
and so on up to
\[ aqr + bpr + cpq = r - 2. \]
If in one of the above steps we find integers $a$, $b$, $c$ such that
\[ aqr + bpr + cpq = i_0, \]
where $1 \leqslant i_0 \leqslant r - 2$, $0 < a < p$, $- q < b < 0$ and $- r <
c < 0$, then we stop the procedure and the triple $\alpha = p - a$, $\beta = -
b$, $\gamma = - c$ solves the Problem with minimum equal to $i_0$. Hence,
$\uplambda^{\tmop{nd}} (f_0) = i_0$ and the deformation $(f_s^{(p - a, - b, -
c)})$ realizes this jump.

If the above search fails, we conclude that $\uplambda^{\tmop{nd}} (f_0) = r -
1$ because the deformation $(f_s^{(b_1, q - a_1, 0)})$, where $a_1 p - b_1 q =
1$, $0 < a_1 < q$, $0<b_1$, realizes this jump.

We may sum up the above considerations in the following
theorem.

\begin{thm}
  \label{twie1}Let $f_0 \in \mathcal{O}_3$ be a convenient and non-degenerate
  singularity with only one two-dimensional face in its Newton diagram. Assume
  that the vertices $(p, 0, 0)$, $(0, q, 0)$, $(0, 0, r)$ of this face are
  such that $p \geqslant q \geqslant r \geqslant 2$ and the numbers $p$, $q$,
  $r$ are pairwise coprime. Then
  \[ \uplambda^{\tmop{nd}} (f_0) = \left\{ \begin{array}{ll}
       i_0 & \begin{array}{l}
         \text{if there exist integers } a, b, c \text{ such that}\\
         aqr + bpr + cpq = i_0, \text{ } 1 \leqslant i_0 \leqslant r - 2,\\
         \text{ } 0 < a < p, \text{ } 0 < - b < q,
         \text{ } 0 < - c < r,\text{ } i_0 \text{ -- minimal},
       \end{array}\\
       & \\
       r - 1 & \text{otherwise}.
     \end{array} \right. \]
  Moreover, $i_0$ can be found algorithmically using only Euclid's algorithm.
\end{thm}

\begin{cor}
  Under the assumptions of Theorem \ref{twie1}, if $r = 2$ then
  $\uplambda^{\tmop{nd}} (f_0) = 1$.
  
\end{cor}

\begin{ex**}
  For $f_0 (x, y, z) \assign x^{11} + y^6 + z^5$ we have $p = 11$, $q = 6$, $r
  = 5$ and{\medskip}
  \begin{center}
  \begin{tabular}{rl}
    $7 \cdot qr - 5 \cdot pr + 1 \cdot pq = 1$ & -- does not satisfy the
    conditions in the theorem\\
    $3 \cdot qr - 4 \cdot pr + 2 \cdot pq = 2$ & -- does not satisfy the
    conditions in the theorem{\hspace{1em}}\\
    $10 \cdot qr - 3 \cdot pr - 2 \cdot pq = 3$ & -- do satisfy the conditions
    in the theorem.
  \end{tabular}
  \end{center}
  
  {\noindent}Hence, $\uplambda^{\tmop{nd}} (f_0) = 3$. This jump is realized by
  the deformation $f_s^{(1, 3, 2)} (x, y, z) \assign x^{11} + y^6 + z^5 + sxy^3
  z^2$. The minimal jump realized by the points lying either in coordinate planes or on axces in equal to $r-1 = 4$.
\end{ex**}

\section*{Aknowledgements}
\noindent Tadeusz Krasinski was supported by OPUS Grant No 2012/07/B/ST1/03293.\\
Szymon Brzostowski and Justyna Walewska were supported by SONATA Grant
NCN No 2013/09/D/ST1/03701.

\bibliographystyle{alpha}

\vspace{5mm}

\noindent {Affiliation/Address}

{\small\it{
  Szymon Brzostowski, Tadeusz Krasi{\'n}ski and Justyna Walewska{
  
  }Faculty of Mathematics and Computer Science{
  
  }University of {\L}{\'o}d{\'z}{
  
  }ul. Banacha 22, 90-238 {\L}{\'o}d{\'z}, Poland{
  
  }brzosts@math.uni.lodz.pl, krasinsk@uni.lodz.pl, walewska@math.uni.lodz.pl
}}

\end{document}